# An Efficient Scenario-based Stochastic Energy Management of Distribution Networks with Distributed Generation, PV Module, and Energy Storage

Mohammad Rasoul Narimani[1], Ali Azizivahed[2], Ehsan Naderi[3]

[1] Electrical and Computer Engineering, Missouri University of Science and Technology, Rolla, USA.
[2] Faculty of Engineering and Information Technology, University of Technology Sydney, Australia.
[3] Electrical Engineering Department, Engineering Faculty, Razi University, Kermanshah, Iran.

**Abstract:** Incorporating Renewable Energy Sources (RES) incurs a high level of uncertainties to electric power systems. This level of uncertainties makes the conventional energy management methods inefficient and jeopardizes the security of distribution systems. In this connection, a scenario-based stochastic programming is introduced to harness uncertainties in the load, electricity price, and photovoltaic generation. Further, a hybrid evolutionary algorithm based on Grey Wolf Optimizer and Particle Swarm Optimisation algorithm is proposed to find the best operation cost, and Energy Not Supplied (ENS) as two important objective functions, which almost always are in stark contrast with each other. The proposed algorithm is applied to the modified IEEE 69-bus test system and the results are validated in terms of efficiency, which indicates a cogent trade-off between the fitness functions addressed above.

## 1. Introduction

Renewable Energy Sources (RESs) have been a key driving force in distribution system analysis in the past few decades. RESs deliver a wide range of benefits to utilities and electricity consumers from economic and sustainability point of views. Electricity consumers can reduce their electricity bills by utilizing solar panels while penetration of these resources can reduce the peak burden on distribution and transmission systems. Beside their benefits, RESs have some drawbacks such as intermittency and fluctuations [1] – [5] which can be alleviated by installing energy storage systems (ESS).

ESSs have recently come under spotlight as a possible mean to reduce the electricity cost, increase the system resiliency, and level out fluctuations in renewable energy resources [6]. Integrating ESSs in distribution power systems as well as accounting RES' intermittencies have forced the conventional analysis techniques to adapt to the new environment. The solution of the deterministic energy management problem cannot be acceptable at the presence of uncertainties incurred to the system by RESs. Therefore, finding new approaches to harness these uncertainties and ameliorating their effect has become an ongoing avenue of research.

Probabilistic analysis is making its way through the management tools, enabling a move from deterministic to stochastic methods, which conduct to robust solutions. One common approach to capture and harness uncertainties in operation and management problem is the scenario based stochastic programming, a promising numerical approach that can handle all source of uncertainty in different problems.

The heart of the scenario based stochastic programming is the probability distribution of forecast errors within which a limited number of scenarios are created to yield the adhered uncertainties of the problem by turning the uncertain problem into an equivalent limited number of deterministic problems. There are multiple ways to create scenarios to solve uncertain problem. Monte Carlo sampling method is a promising approach that creates scenarios through the roulette wheel mechanism which is based on the probability Density functions (PDF) of random input variables. Higher number of scenarios can capture uncertainties more accurately but incurs higher computational burden to the problem. Thus, there should be a compromise between precision and computational burden which in this paper is observed by a scenario reduction method and a criterion which decide whether the current number of scenarios yield pre-specified level of precision.

Distribution System Operator (DSO) used to consider the operation cost as their priority. However, reliable electricity is as fundamental as the affordable electricity in modern society in which pushes the system operator to operate the system not only at the most economical way but also at a reasonable level of reliability. Reliability measures by frequency, duration, and magnitude of adverse effects on the system [7]. In this connection, Energy Not Supplied (ENS)—a prominent reliability index based on frequency, duration, and magnitude of adverse effects–is considered as an objective function beside the operation cost at this paper which turns the problem to a Multi-Objective Optimisation Problem (MOOP). A multi-objective method and decision making strategy have been employed to handle the MOOP and find the Best Compromise Solution (BCS), respectively.

A plethora of researches has been proposed to solve the energy management problem in distribution systems from different aspects [8]-[13]. Farzin et al. proposed a Mixed Integer Linear Programming (MILP) for solving stochastic energy management problem of micro-grids during islanding events [14]. Eajal et al. presented a stochastic energy management and Unit Commitment (UC) in AC/DC smart grid considering Plug-in Electric Vehicles (PEVs), distributed generation sources, and BEESs [15]. Li et al. proposed a Stochastic Dynamic Programming (SDP) to solve an optimal energy management problem for Plug-in Hybrid Electric Bus (PHEB) [16]. Wang et al. presented a stochastic approach for energy management in distribution networks considering DGUs and a feed-in tariff program

[17]. Su et al. introduced a stochastic programming for energy scheduling in micro-grids considering RERs and PEVs [18]. A SDP approach is presented in [19] by Wu et al. for smart home energy management at the presence of PEVs, PV arrays, and BESSs.

Considering the reliability factor and solving the problem as a MOOP as well as employing scenario based approach for solving the proposed problem distinguishes the proposed study from those in literature. Last but not least, solving the stochastic-based energy management problem as a MOOP needs a powerful and efficient optimisation algorithm. This need is satisfies by developing a novel Hybrid Grey-Wolf Optimizer-Particle Swarm Optimisation (HGWO-PSO) algorithm. Yet another, the proposed problem is a MOOP which has a set of optimal solutions (Pareto-optimal fronts) instead of one unique solution. A Fuzzy decision making engine is employed to select the best compromise solution based on the operator desires.

The contribution of this paper is three fold, which the qualities itemized hereunder:
- Presenting a stochastic framework for the energy management problem.
- Taking reliability objective into consideration and solving the problem as a MOOP in stochastic environment.
- Proposing a new hybrid evolutionary algorithm.

## 2. Problem Formulation

A Brief description about the proposed stochastic framework accompany with detailed formulation of the proposed Multi-Objective Energy Management (MOEM) problem are presented in this section.

### 2.1. Scenario-Based Stochastic Framework

There are many numbers of distributions in statistical analysis, normal distribution is the most common one to handling uncertainties in stochastic studies. In this study, normal distribution is chosen to reverberate the strong influences of uncertainties of Renewable Energy Sources (RESs).

Additionally, it is common to generate different number of scenarios to precisely handle problems' uncertainties. It seems axiomatic that more scenarios beats fewer numbers so much so that they can definitely enhance accuracy of the process. Yet this strategy will pose too much CPU time for computing; therefore, not only is it necessary to execute the problem for only a limited number of scenarios, but also it is crucial to rule out the extra ones in an expeditious way. However, interested readers are directed to [20], [21] for more information about aforementioned processes.

To sum up, in the subsequent sections, the probability of $s^{th}$ scenario at $t^{th}$ time interval is marked as $\psi_s$. Furthermore, numerical value of loads, power generated of PV arrays, and price of electrical energy of network constitute the total number of uncertain variables in this study.

### 2.2. Energy Management Formulation

In this section objective functions and all constraints and limitations are elaborated upon hereunder.

*2.2.1. Objective Function 1 (Total Operation Cost):* The total operation cost in distribution systems can be captured by (1). Let $\psi_s$ denotes the probability of $s^{th}$ scenario at $t^{th}$ time interval.

$$F_1 = \sum_{s=1}^{N_s} \psi_s \times Cost_s$$
$$= \sum_{s=1}^{N_s} \psi_s \left( \sum_{i=1}^{N_{SS}} \sum_{t=1}^{T} EP_{SS_{i,s}}^t \times RE_{SS_{i,s}}^t \right. \quad (1)$$
$$\left. + \sum_{j=1}^{N_{RES}} \sum_{t=1}^{T} EP_{RES_{j,s}}^t \times RE_{RES_{j,s}}^t \right)$$

where, $RE_{SS_{i,s}}^t$ and $EP_{SS_{i,s}}^t$ are the receiving energy and energy price from/at $i^{th}$ sub-station in scenario $s$ at $t^{th}$ time interval, respectively; $RE_{RES_{j,s}}^t$ and $EP_{RES_{j,s}}^t$ are the receiving energy and energy price from/of $j^{th}$ RES in scenario $s$ at $t^{th}$ time interval, respectively; $N_s$ and $N_{SS}$ are the number of scenarios and sub-stations; $T$ is the total number of time intervals in one day i.e. 24 hours; $N_{RES}$ is the total number of RESs; and $F_1$ is the total operation cost of distribution system ($).

*2.2.2. Objective Function 2 (ENS):* There are different indices to evaluate reliability of system, which are usually defined by frequency and duration of interruption events. Contrary to most studies, a reliability index i.e. ENS, which takes both frequency and duration of interruption into account, is considered in this paper which can be formulated as (2).

$$F_2 = \sum_{s=1}^{N_s} \psi_s \times ENS_s$$
$$= \sum_{s=1}^{N_s} \psi_s \left( \sum_{b=1}^{N_B} P_{b,s} \right. \quad (2)$$
$$\left. \times \sum_{r \in R} AT_{r,b,s}^{repair} + AT_{r,b,s}^{restoration} \right)$$

where, $AT_{r,b,s}^{repair}$ and $AT_{r,b,s}^{restoration}$ are the annual repair and restoration times related to all branches associated with bus $b$ in scenario $s$, respectively; $R$ includes all buses that feed from specific feeder; $N_B$ is the number of buses in the system; and $P_{b,s}$ is the active power of $b^{th}$ bus in scenario $s$.

*2.2.3. Cost-Profit Analysis:* A cost-profit analysis is conducted to recognize when the investment cost will be returned or compensated. The cost-benefit analysis for implemented PV arrays accompany with their build-in ESSs can be considered as (3) & (4).

$$Profit = C_{NPV} \sum_{s=1}^{N_s} \psi_s \times \sum_{t=1}^{365} \Delta(TOC_s) \quad (3)$$

$$\Delta(TOC_s) = TOC_s^{old} - TOC_s^{new} \quad (4)$$

where, $TOC_{new}^s$ and $TOC_{old}^s$ are the total operation cost with and without considering PV cells and ESSs in scenario s; $C_{NPV}$ is a coefficient called net present value, equals to 1.07, and defined by the difference between the present value of cash inflows and the present value of cash outflows over a period of time; and $Profit$ indicates the obtained profit.

*2.2.4. Constraints and Restrictions:* All constraints in the proposed energy management problem are represented as (5)-(13).

$$\forall s \in N_s \Rightarrow P_{SS}^{t,s} + P_{PV}^{t,s} + P_{ESS}^{t,s} + P_{DG}^{t,s} = P_{Loss}^{t,s} + P_{Demand}^{t,s} \quad (5)$$

$$\forall s \in N_s, \forall t \in T, \forall k \in ESS \Rightarrow W_{k,t,s}^{ESS} = W_{k,t-1,s}^{ESS} + \left(\xi_{C,k} \times P_{C,k,t,s} \times \Delta t\right) - \left(\frac{1}{\xi_{D,k}} \times P_{D,k,t,s} \times \Delta t\right) \quad (6)$$

$$\forall s \in N_s, \forall t \in T \Rightarrow P_{m,t,s} = \sum_{m=1}^{N_B} (V_{m,t,s} \times V_{n,t,s} \times Y_{mn}) \times \cos(\theta_{mn,s} - \delta_{m,t,s} - \delta_{n,t,s}) \quad (7)$$

$$\forall s \in N_s, \forall t \in T \Rightarrow Q_{m,t,s} = \sum_{m=1}^{N_B} (V_{m,t,s} \times V_{n,t,s} \times Y_{mn}) \times \sin(\theta_{mn,s} - \delta_{m,t,s} - \delta_{n,t,s}) \quad (8)$$

$$\forall s \in N_s, \forall t \in T, \forall k \in N_{ESS} \Rightarrow W_{k,min}^{ESS} \leq W_{k,t,s}^{ESS} \leq W_{k,max}^{ESS} \quad (9)$$

$$\forall s \in N_s, \forall t \in T, \forall k \in N_{ESS} \Rightarrow P_{C,k,t,s} \leq P_{C,k}^{max} \quad (10)$$

$$\forall s \in N_s, \forall t \in T, \forall k \in N_{ESS} \Rightarrow P_{D,k,t,s} \leq P_{D,k}^{max} \quad (11)$$

$$\forall s \in N_s \Rightarrow S_{l,s} \leq S_l^{max} \quad (12)$$

$$\forall s \in N_s, \forall t \in T \Rightarrow V_m^{min} \leq V_{m,t,s} \leq V_m^{max} \quad (13)$$

where, $s$, $t$, and $k$ are three indices, which are related to number of scenarios, time intervals, and number of ESSs, respectively; $N_s$ is the number of scenarios; $T$ is the total number of time intervals in one day–24 hours; $N_{ESS}$ is the total number of employed ESSs; $P_{SS}^{t,s}$ is the power that receives from sub-station in $s^{th}$ scenario at $t^{th}$ time interval; $P_{PV}^{t,s}$ is the power of PV array that injects into the system in $s^{th}$ scenario at $t^{th}$ time interval; $P_{ESS}^{t,s}$ is the power of ESS in $s^{th}$ scenario at $t^{th}$ time interval; $P_{DG}^{t,s}$ is the power of diesel generator that injects into system in $s^{th}$ scenario at $t^{th}$ time interval; $P_{Loss}^{t,s}$ and $P_{Demand}^{t,s}$ are power losses and power demand of the distribution system in $s^{th}$ scenario at $t^{th}$ time interval; $W_{k,t,s}^{ESS}$ and $W_{k,t-1,s}^{ESS}$ are the levels of energy within the $k^{th}$ ESS in $s^{th}$ scenario at $t^{th}$ and $(t-1)^{th}$ time intervals, respectively; $P_{C,k,t,s}$ and $P_{D,k,t,s}$ are permissible rates of charge and discharge of $k^{th}$ ESS during a specified time period, one hour, in $s^{th}$ scenario at $t^{th}$ time interval, respectively; $\xi_{C,k}$ and $\xi_{D,k}$ are the efficiency of the $k^{th}$ ESS during charge and discharge processes, respectively; $P_{m,t,s}$ and $Q_{m,t,s}$ are active and reactive power that inject by $m^{th}$ bus in $s^{th}$ scenario at $t^{th}$ time interval, respectively; $V_{m,t,s}$ and $V_{n,t,s}$ are the voltage amplitudes of buses m and n in $s^{th}$ scenario at $t^{th}$ time interval, respectively; $\delta_{m,t,s}$ and $\delta_{n,t,s}$ are phase angles of voltage of buses m and n in $s^{th}$ scenario at $t^{th}$ time interval, respectively; $Y_{mn}$ is the magnitude of admittance between buses m and n; $\theta_{mn,s}$ is phase angle of admittance between buses m and n in $s^{th}$ scenario; $S_{l,t,s}$ and $S_l^{max}$ are the amount of power flow and its corresponding maximum bounds related to $l^{th}$ branch in $s^{th}$ scenario at $t^{th}$ time interval, respectively; $W_{k,min}^{ESS}$ and $W_{k,max}^{ESS}$ are the minimum and maximum amounts of energy that can be stored in $k^{th}$ ESS, respectively; $P_{C,k}^{max}$ and $P_{D,k}^{max}$ are the maximum rate of charge and discharge processes of $k^{th}$ ESS during each time interval, respectively; and $V_m^{min}$ and $V_m^{max}$ are minimum and maximum bounds of voltage magnitude, respectively.

## 3. Multi-objective Solution Methodology

In this section only a brief explanations are provided about the proposed multi-objective methodology and its essential processes. Interested readers are directed to [22]-[33] for more information.

### *3.1. Trapezoidal Membership Function*

A trapezoidal membership function is employed in order to normalize the numeral value of each objective function, which can be written as (14).

$$\forall h \in \{1,2\} \Rightarrow \Psi_h(X) = \begin{cases} 0 & F_h \geq F_h^{min} \\ \frac{F_h^{max} - F_h}{F_h^{max} - F_h^{min}} & F_h^{min} \leq F_h \leq F_h^{max} \\ 1 & F_h \leq F_h^{max} \end{cases} \quad (14)$$

where, $X$ is the vector of decision variables; $F_h$ is $h^{th}$ objective function, $h = 1,2$; $F_h^{min}$ and $F_h^{max}$ are minimum and maximum bounds of $h^{th}$ objective function, respectively; and $\Psi_h$ is the fuzzy set for $h^{th}$ objective function.

### *3.2. Non-dominance Concept and Pareto-technique*

In contrast to the single-objective optimisation problem which always has a unique optimal solution, MOOPs have a set of non-dominated optimal solutions, Pareto-optimal front. In this connection, the dominance concept, based on (15) & (16), is employed to find all the optimal solutions. Note that, $X$ is the vector of decision variables.

$$\forall x = \{1,2\} \Rightarrow F_x(X_1) \leq F_x(X_2) \quad (15)$$
$$\exists y = \{1,2\} \Rightarrow F_y(X_1) < F_y(X_2) \quad (16)$$

### *3.3. Best Compromise Solution (BCS)*

Equation (17) is applied to all the retained optimal solutions to select the Best Compromise Solution (BCS);

$$\Upsilon_\Psi = \frac{\sum_{h=1}^2 \{\Omega_h \times \Psi_{qh}\}}{\sum_{q=1}^{N_{POF}} \sum_{h=1}^2 \{\Omega_h \times \Psi_{qh}\}} \quad (17)$$

where, $N_{POF}$ is the number of non-dominated solutions; $\Omega_h$ is the weight factor for $h^{th}$ objective function; $\Psi_{qh}$ is the corresponding normalized objective function; and $\Upsilon_\Psi$ is the BCS.

## 4. Proposed Optimisation Algorithm

This section provides a brief explanation about the proposed hybrid algorithm.

### 4.1. Overview of Grey Wolf Optimizer

GWO mechanism is inspired by the hunting process of grey wolves [34]. First, a population of grey wolves is generated randomly, which are a set of candidate solutions for the current optimisation problem. Leader of the wolves, the individual with best fitness function, is called $\alpha$. The next wolf in terms of priority, which is the leader's advisor is called $\beta$, which can help the leader in decision-making processes. Wolves with low ranks of fitness are called $\omega$. $\omega$ wolves with better positions are called $\delta$ and exchange information with both $\alpha$ and $\beta$. Mathematical model for updating the position of low ranks solutions, $\omega$ wolves with low fitness, can be represented as (18) & (19).

$$Dis = |\{\eta \otimes X_{prey}^{old}\} - X_{wolf}^{old}| \tag{18}$$
$$X_{wolf}^{New} = X_{prey}^{old} - \{\zeta \otimes Dis\} \tag{19}$$

where, $X_{wolf}^{old}$ indicates the vector of wolves' position in current iteration; $X_{wolf}^{New}$ indicates the vector of wolves' position in next iteration; $X_{prey}^{old}$ indicates the position vector of prey; $\zeta = \varepsilon(2R_2 - 1)$; $\eta = 2R_1$; $\varepsilon$ is linearly decreased from two to zero over the course of iterations; and $R_1$ and $R_2$ are random numbers in the range of zero to one.

$\omega$ wolves modify and update their situations with regard to $\alpha$, $\beta$, and $\delta$ according to equations (20)-(26) in which $Dis_\alpha$, $Dis_\beta$, and $Dis_\delta$ mark the distance between hunt and $\alpha$, $\beta$, and $\delta$, respectively.

$$Dis_\alpha = |\{\eta_1 \otimes X_\alpha^{old}\} - X_{wolf}^{old}| \tag{20}$$
$$Dis_\beta = |\{\eta_2 \otimes X_\beta^{old}\} - X_{wolf}^{old}| \tag{21}$$
$$Dis_\delta = |\{\eta_3 \otimes X_\delta^{old}\} - X_{wolf}^{old}| \tag{22}$$
$$X_{wolf,1} = X_\alpha^{old} - \{\zeta_1 \otimes Dis_\alpha\} \tag{23}$$
$$X_{wolf,2} = X_\beta^{old} - \{\zeta_2 \otimes Dis_\beta\} \tag{24}$$
$$X_{wolf,3} = X_\delta^{old} - \{\zeta_3 \otimes Dis_\delta\} \tag{25}$$
$$X_{wolf}^{New} = \frac{\{X_{wolf,1} + X_{wolf,2} + X_{wolf,2}\}}{3} \tag{26}$$

where, $\eta_1 - \eta_3$ and also $\zeta_1 - \zeta_3$ are coefficients vectors that generate randomly.

According to [34], GWO can precisely carry out both exploration and exploitation processes in the search space by adjusting $\zeta$ & $\eta$ parameters. In half of iterations while $|\zeta| > 1$ the GWO is exploring and in the other half when $|\zeta| < 1$ the algorithm is exploiting.

### 4.2. Overview of Particle Swarm Optimisation

Like every other evolutionary algorithm, in the PSO at first a population consist of $Size_{Pop}$ members is generated randomly. Each member in the population promotes toward the global optimal solution by adopting its features–position and velocity–in the search space. In this mobility, position and velocity of each individual can be updated using (27) & (28).

$$\forall i \in \{1,2,\dots,Size_{Pop}\} \Rightarrow v_i^{New} = v_i^{old} + v_i^{New} \tag{27}$$

$$v_i^{New} = \mu v_i^{old} + c_1 R_1 (\varphi_i^{old} - v_i^{old}) + c_2 R_2 (\Phi^{old} - v_i^{old}) \tag{28}$$

where, $v_i^{old}$ and $v_i^{New}$ indicate current and updated position of $i^{th}$ particle, respectively; $v_i^{old}$ and $v_i^{New}$ indicate current and updated velocity of $i^{th}$ particle, respectively; $\varphi_i^{old}$ is the best fitness by the objective function at current iteration; $\Phi^{old}$ is the global best position at current iteration; and $\mu$ is the inertia coefficient of particles.

Interested readers are referred to [35] for getting detailed information about different models of PSO algorithm.

### 4.3. Proposed Hybrid GWO-PSO Approach

For brevity, the hybridization process outlined hereunder;
- The initial population splits into two equal sections in terms of individuals, each of which includes half of the population.
- Each section promotes through one of GWO or PSO algorithms. In this connection, each algorithm handles the constraints independently.
- The best solution is obtained at each iteration, and the population is shuffled together.
- The best solution for both GWO and PSO algorithms are updated and fed into both algorithms.

The entire population splits into two sets randomly, and this process will be repeated until the algorithm converges to the global solution.

## 5. Simulation and Numerical Results

This section interprets the obtained numerical results for solving both deterministic (Case I) and stochastic (Case II) approaches in handling the adhered uncertainties in the proposed problem. Moreover, four various sizes of scenarios are considered in the stochastic optimization in order to demonstrate the effect of size of scenarios on the expected value. It should be noted that each scenario consists of the output power for PV arrays, load of system, and energy price in 24-hours.

### 5.1. Initialization

The proposed algorithm has several parameters, which must be tuned in advance for yielding balance between both exploration and exploitation processes.

The population size is set to 500; the maximum number of iteration is set to 50; the interval of inertia weight is bounded into [0.4, 0.9]; and the values of both learning factors are set to 1.49618.

The whole scheme for solving the proposed problem is implemented in MATLAB R2014 environment using quad-core processor laptop machine with 1.6 GHz clock frequency and 4.0 GB of RAM.

The proposed algorithm is applied on IEEE 69-bus distribution system (see Fig. 1) which is retrofitted by diesel generators, PV panels, and ESSs. Detailed data for this test system can be found in [36].

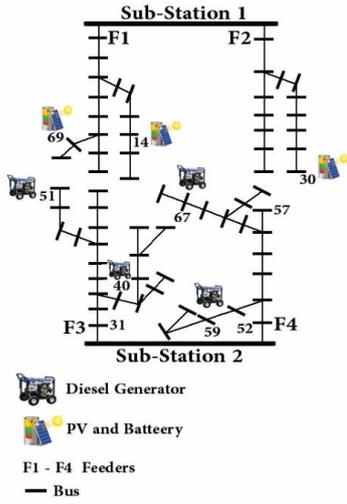

**Fig. 1.** *Single-line diagram of modified IEEE 69-bus distribution network with DGs, PVs, and ESSs*

According to Fig. 1, three units of PV arrays with capacity of 1.5 MW accompanied with their ESSs are placed at buses #14, #30 and #69. In addition, DGs are installed at buses #40, #51, #59 and #67. The capital expense of PV arrays, converters as well as ESSs are equal to $2,000/KW, $400, and $100/KW, respectively. The total investment cost of both PV arrays and ESSs during 20-year period (ESS's lifetime is five years and they should be replaced four times in the proposed time horizon) is $9,751,200. Note that, the horizon of this study is a daylong with hourly steps.

Both versions of study, Cases I and II, are elaborated upon in the following section so as to scrutinize the ability of the proposed methodology in handling uncertainties in which the effect of scenario's size on the expected value is highlighted.

### 5.2. Obtained Results and Analysis

The best obtained results in solving single- and multi-objective version of the problem are tabulated in Table 1 for Case I. Note that, the obtained results in multi-objective approach are the Best Compromise Solution (BCS).

**Table 1** Best single- and multi-objective results for Case I

| Models | Objectives | Operation Cost ($) | ENS (KWh/year) |
|---|---|---|---|
| SO[1] | Cost | 2727.94 | 2067224.77 |
|  | ENS | 2977.79 | 1546895.58 |
| MO[2] | BCS | 2828.49 | 1864182.51 |

[1]SO, Single-objective; [2]MO, Multi-objective

According to Table 1 it is clear that, the proposed algorithm could compromise between different objective functions since the BCS lays between the best optimal solutions obtained in single-objective optimisation problems for different objective functions. Unfortunately, there was no literature to compare the obtained results with. Since the prominent contribution of this paper stems from handling and elaborating uncertainties' effects on MOEM problem, the rest of this paper covers stochastic related analysis.

The best obtained results in solving single- and multi-objective version of the problem for Case II are tabulated in Tables 2 & 3, respectively. Different well-known statistical factors including average values, Standard Deviation (SD), 95% Confidence Intervals (CIs), Expected Values (EVs), and Relative Errors (REs) are employed to analyze the obtained results. Comparing the objective function values in Tables 2 & 3 demonstrates that the proposed approach can handle the problem even in the stochastic environment.

**Table 2** Best single-objective obtained results for Case II

| Obj. | Scenarios | Mean | SD | 0.95 CI | EV |
|---|---|---|---|---|---|
| Cost | 1 × 30 | 2747.70 | 11.06 | 3.96 | 2740.51 |
|  | 2 × 30 | 2748.90 | 11.92 | 3.69 | 2749.54 |
|  | 3 × 30 | 2730.50 | 11.77 | 2.43 | 2748.01 |
|  | 4 × 30 | 2717.30 | 11.89 | 2.13 | 2731.90 |
| ENS | 1 × 30 | 1570200.00 | 7863.20 | 2813.80 | 1544247.36 |
|  | 2 × 30 | 1545000.00 | 7995.30 | 2477.80 | 1543119.92 |
|  | 3 × 30 | 1553512.88 | 8306.42 | 1716.12 | 1548153.93 |
|  | 4 × 30 | 1537995.50 | 8261.51 | 1478.17 | 1539900.00 |

**Table 3** Best multi-objective results (BCSs) for Case II

| Obj. | Scenarios | Mean | SD | 0.95 CI | EV |
|---|---|---|---|---|---|
| Cost | 1 × 30 | 2836.30 | 10.54 | 3.77 | 2848.69 |
|  | 2 × 30 | 2834.30 | 11.45 | 3.55 | 2840.38 |
|  | 3 × 30 | 2826.80 | 11.12 | 2.30 | 2825.65 |
|  | 4 × 30 | 2847.20 | 11.12 | 2.05 | 2834.33 |
| ENS | 1 × 30 | 1858100.00 | 7863.20 | 2813.80 | 1864818.56 |
|  | 2 × 30 | 1858300.00 | 7995.30 | 2477.80 | 1867638.92 |
|  | 3 × 30 | 1880950.00 | 8306.42 | 1716.10 | 1873753.81 |
|  | 4 × 30 | 1877325.00 | 8497.18 | 1595.20 | 1871006.79 |

The mean value of objective function in Tables 2 & 3 is decreased by increasing the number of scenarios. This trend means that the algorithm has converged to lower values of cost and ENS by increasing the number of scenarios. Conversely, SD value that measures the spread of samples is increased as it was expected. The nature of scenario generation is completely randomly. Thus considering more scenarios implies taking more outliers into account, and increase the SD value.

Confidence interval, a range computed using sample statistics to estimate an unknown population parameter with a stated level of confidence, plays a crucial role in statistical analysis. The desired level of confidence is set to 95%. The lower value of confidence interval demonstrates a better precision for the expected value. Increasing number of scenarios in Tables 2 & 3 decreases the value of 0.95 CI but incurs computational burden to the problem. Therefore, leveraging a tradeoff between precision and computational is in most of interest which is yield in this paper by employing the stopping rule in [21].

Comparing results in the Tables 1 & 2 demonstrates that, taking uncertainties into account increases the total cost and ENS objective functions compared to their corresponding values in deterministic approach. It was expected because there are different scenarios in Case II, instead of one fixed scenario in the Case I, which can potentially cause deteriorating the optimal values.

The optimal charging and discharging patterns for ESSs in 24-hour time window are depicted in Fig. 2 for different objective functions in Cases I & II. Positive and negative values represent charging and discharging cycles, respectively. From Fig. 2 it is clear that the ESSs charge during lightly-loaded hours, where the electricity price is low, and discharge during the heavily-loaded hours, where the electricity price is high. These charging and discharging patterns helps distribution network to lighten the peak stress.

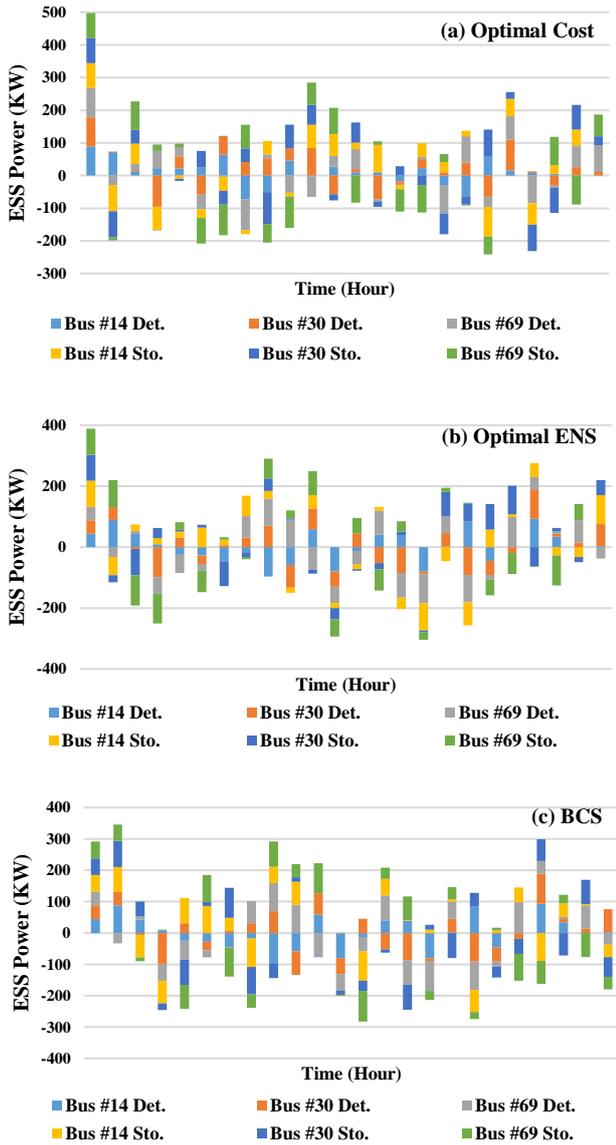

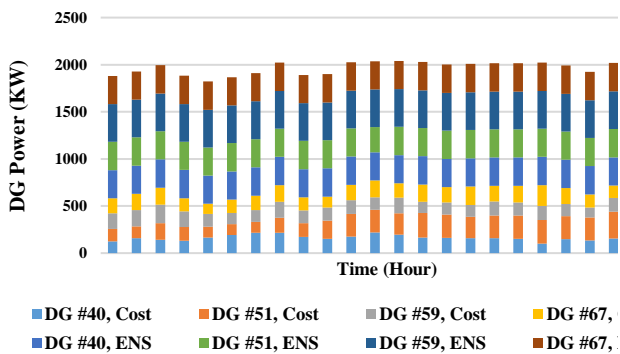

**Fig. 2.** *Obtained power for all ESSs related to Cases I & II* (**a**) *Optimal cost,* (**b**) *Optimal ENS, and* (**c**) *Best compromise solution*

The optimal level of power generated of DGs for cost-based and ENS-based energy management problems, Cases I & II, are depicted in Figs. 3 & 4, respectively. Furthermore, the output power of DG units related to the BCSs, when both Cases consider simultaneously, are depicted in Fig. 5.

**Fig. 3.** *Obtained active power of DG units for Case I*

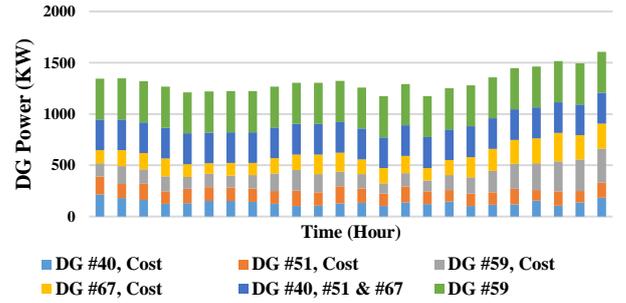

**Fig. 4.** *Obtained active power of DG units for Case II and scenario 4 × 30*

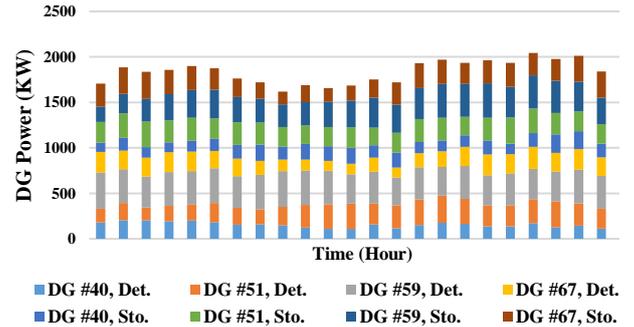

**Fig. 5.** *BCS of DG units for Case I & Case II (scenario 4 × 30)*

According to Figs. 3-5, DG outputs are different for different objective function as it was expected. The generation levels of DGs are low in cost-based energy management problem because it is more economical to drawn power from the grid rather than generating by DGs. Conversely for the reliability-based energy management DG outputs are high because distributed generation can improve the ENS index. Moreover, the DG outputs in Fig. 5 lay between their corresponding values in Figs. 3 & 4 for both deterministic and stochastic versions as it was expected. Yet another, difference between DG outputs in Cases I & II underscores the importance of considering uncertainties in energy management problem. Because the optimal decision variables for deterministic problem is not necessarily the optimal solution in stochastic environment. The drawn power from sub-station has been illustrated in Fig. 6 for all objective and case studies.

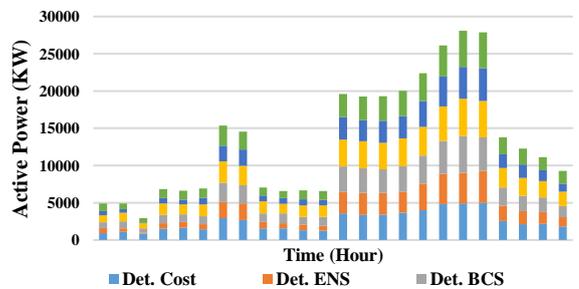

**Fig. 6.** *Plots of sub-station active power in different aspects*

Although there is a slight difference between drawing powers from substation for the same objective in Case I and Case II in Fig. 6, as it was an expected they follow a same pattern.

The Pareto-Optimal Fronts (POFs) are depicted in Fig. 7 for both deterministic and stochastic cases. Despite the complexity of the problem the proposed algorithm could obtain a well-distributed POF in both cases. It is notable that, the number of non-dominated solutions has been decreased in stochastic case as it was expected. This stems from eliminating the low probable points in the search space of the problem. Moreover, all the obtained points in POFs are non-dominated solutions, which prove the suitability of the proposed algorithm for handling the complex MOOPs in both deterministic and stochastic environment.

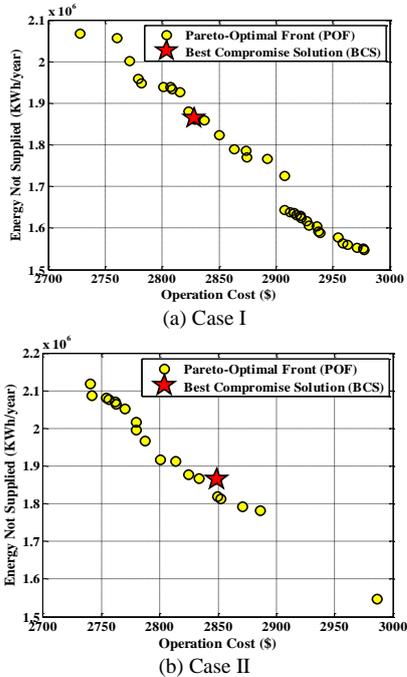

**Fig. 7.** *Two-dimensional POFs from different aspects* **(a)** *Deterministic solution,* **(b)** *Stochastic solution related to scenario* $1 \times 30$

The Probability Density Functions (PDFs) for optimal operation cost and ENS objective function are depicted in Fig. 8 for stochastic single-objective energy management problem.

It is expected that objective functions follow the patterns of input random variables. In this connection a normal distribution is mapped into both pictures in Fig. 8 in order to show how close the obtained patters are to the pattern of the input random variables. In other words, the output and inputs PDFs are similar, which means that the proposed method can handle uncertainties without changing inherit of the problem.

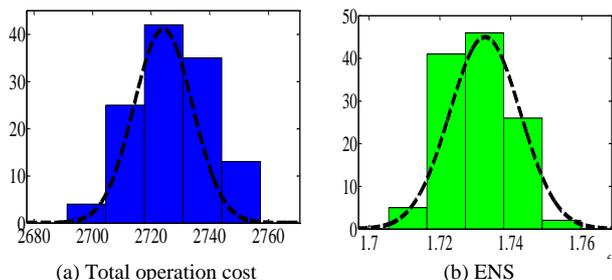

**Fig. 8.** *PDF for optimizing total operation cost & ENS objective functions*

Moreover, PDFs corresponds to the BCS are illustrated in Fig. 9 for both objectives. From this figure it is clear the objective functions follow the patterns of input random variables even in multi-objective environment. It also indicates that the proposed algorithm handle MOOP in stochastic environment without changing the problem's characteristics.

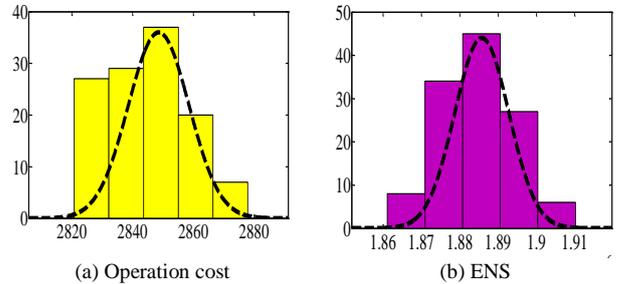

**Fig. 9.** *PDF related to BCS* **(a)** Operation cost, and **(b)** ENS

Financial justification in every engineering project plays a vital role and is in most of interest. In this connection the return of investment of the proposed approach is analyzed. Fig. 10 shows the cumulative profit value for a 20-year time horizon. It is crystal clear that, the investment cost will return after 12 years and the net profit is equal to $12,566,086.

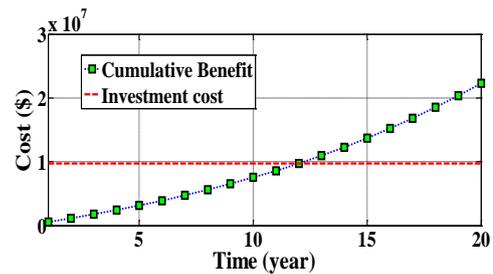

**Fig. 10.** *The profit and investment costs for a 20-year time horizon*

## 6. Conclusions

An efficient hybrid optimisation algorithm has been successfully implemented for solving single- and multi-objective deterministic and stochastic energy management problems in distribution systems. In this connection, the proposed scenario-based approach is implemented to harness the adhered uncertainties to the problem. A tradeoff between precision and computational burden is yielded. The proposed hybrid evolutionary algorithm is applied on IEEE 69-bus test system and the obtained results prove the algorithms' ability in solving MOEM problem in both deterministic and stochastic environments. Furthermore, the proposed cost-profit analysis adds more value to this study and makes it conspicuous from other literature in the area.